\documentclass[11pt]{article}

\usepackage[utf8]{inputenc}
\usepackage{amsmath}
\usepackage{amsfonts,amssymb,latexsym,multicol,color,mathtools}
\usepackage[english]{babel}
\usepackage[T1]{fontenc}

\newcounter{fig}
\newtheorem{prop}{Proposition}

\newcommand{\expli}[1]{\quad\text{\footnotesize (#1)}}
\newcommand{\epv}{\quad ; \quad}

\makeatletter
\ifcase\@ptsize

\or

\or

\fi
\makeatother

\newcommand{\abs}[1]{\left\lvert #1 \right\rvert}

\newcommand{\set}[1]{\{ #1 \}}

\newcommand{\eps}{\varepsilon}

\newcommand{\ioe}{\leqslant}
\newcommand{\soe}{\geqslant}

\newcommand{\vers}{\rightarrow}

\newcommand{\demi}{{\frac{1}{2}}}

\newcommand{\Dcal}{{\mathcal D}}
\newcommand{\Ecal}{{\mathcal E}}

\newcommand{\Pcal}{{\mathcal P}}

\newcommand{\norm}[1]{\left\| #1 \right\|}

\newcommand{\Nat}{{\mathbb N}}

\newcommand{\Real}{{\mathbb R}}

\newcommand{\virg}{\raisebox{.7mm}{,}}

\newcommand{\fin}{\hfill$\Box$}
\newcommand{\dem}{\noindent {\bf Proof\ }}
\newcommand{\fine}{\tag*{\mbox{$\Box$}}}

\title{On the infimum of the absolute value\\ of successive derivatives\\ of a real function defined on a bounded interval}

\author{Michel Balazard}

\begin{document}
\maketitle

\begin{center}
  {\sc Abstract}
\end{center}
\begin{quote}
{\footnotesize 
A study of the greatest possible ratio of the smallest absolute value of a higher derivative of some function, defined on a bounded interval, to the $L^p$-norm of the function.}
\end{quote}

\begin{center}
  {\sc Keywords}
\end{center}
\begin{quote}
{\footnotesize Chebyshev polynomials, Legendre polynomials, extremal problems, inequalities for derivatives  }
\end{quote}



\begin{flushright}
To the memory of Eduard Wirsing,\\ master of analysis,\\ and of its applications to number theory.
\end{flushright}

\section{Introduction}

Let $n$ be a positive integer, $I=[a,b]$ a bounded segment of the real line, of length $L=b-a$. Define $\Dcal^n(I)$ as the set of real functions $f$ defined on $I$, with successive derivatives $f^{(k)}$ defined and continuous on $I$ for $0 \ioe k \ioe n-1$, and $f^{(n)}$ defined on $\mathring{I}=\,]a,b[$. We will use the notation
\[
m_n(f)=\inf_{a<t<b} \vert f^{(n)}(t)\vert.
\]

Let $p$ be a positive real number, or $\infty$.

The problem addressed in this article is that of determining the best constant $C^*=C^*(n,p,I)$ in the inequality
\[
m_n(f) \ioe C^*\norm{f}_p \quad (f \in \Dcal^n(I)),
\]
where
\[
\norm{f}_p =\Big(\int_a^b \abs{f(t)} ^p \, dt \Big)^{1/p},
\]
with the usual convention when $p=\infty$, here: $\norm{f}_{\infty}=\max \abs{f}$.

\smallskip

This problem has been posed by Kwong and Zettl in their 1992 Lecture Notes \cite{zbMATH00205893} (see Lemma~1.1, p. 6). They give upper bounds for $C^*(n,p,I)$, but their reasoning and results are erroneous. In her 1993 PhD Thesis \cite{MR2689365}, Huang has pointed out that this problem is equivalent to a classical problem in the theory of polynomial approximation: that of determining the minimal~$L^p$-norm of a monic polynomial of given degree on a given bounded interval. Our purpose in this text is to give a     new proof of the equivalence, and  to list the consequences of the known results about this extremal problem for the evaluation of~$C^*(n,p,I)$.

\section{First observations}

\subsection{Homogeneity}

Defining $g(u)=f(a+uL)$ for $f \in \Dcal^n(I)$ and $0\ioe u \ioe 1$, one has
\[
g \in \Dcal^n([0,1]) \epv g^{(n)}(u)=L^n f^{(n)}(a+uL) \quad (0 <u<1) \epv \norm{g}_p=L^{-1/p}\norm{f}_p.
\]
Hence,
\begin{equation}\label{220817b}
C^*(n,p,I)=C^*(n,p,[0,1])\,L^{-n-1/p},
\end{equation}
and one is left with determining $C^*(n,p,[0,1])=C(n,p)$, or in fact $C^*(n,p,I)$ for any fixed, chosen segment $I$. We will see that $I=[-1,1]$ is particularly convenient.

\subsection{An extremal problem}

One has
\begin{align*}
C^*(n,p,I)&=\sup\set{m_n(f)/\norm{f}_p, \; f\in \Dcal^n(I), \; m_n(f) \neq 0}\\
&=\sup\set{m_n(f)/\norm{f}_p, \; f\in \Dcal^n(I), \; m_n(f)= \lambda}\expli{for every $\lambda >0$}\\
&=\lambda/D^*(n,p,\lambda, I),
\end{align*}
where
\begin{align*}
D^*(n,p,\lambda, I)&=\inf\set{\norm{f}_p, \; f\in \Dcal^n(I), \; m_n(f)=\lambda}\\
&=\inf\set{\norm{f}_p, \; f\in \Dcal^n(I), \; m_n(f)\soe\lambda},
\end{align*}
the last equality being true since $D^*(n,p,\mu, I)=\frac{\mu}{\lambda} D^*(n,p,\lambda, I)\soe D^*(n,p,\lambda, I)$ if $\mu \soe \lambda$. 
\smallskip

Also, since a derivative has the intermediate value property (cf. \cite{zbMATH02715424}, pp. 109-110), the inequality $m_n(f)\soe \lambda>0$ implies that $f^{(n)}$ has constant sign on $I$, so that
\begin{equation*}
D^*(n,p,\lambda, I)=\inf\set{\norm{f}_p, \; f\in \Dcal^n(I), \; f^{(n)}(t) \soe\lambda \text{ for } a<t<b}.
\end{equation*}

Thus, determining $C^*(n,p,I)$ is equivalent to minimizing $\norm{f}_p$ for $f \in \Dcal^n(I)$ with the constraint $f^{(n)}(t) \soe\lambda>0$ for $a<t<b$. We will denote this extremal problem by $\Ecal^*(n,p,\lambda, I)$.

\section{The relevance of monic polynomials}

Let $\Pcal_n$ be the set of monic polynomials of degree $n$, with real coefficients, identified with the set of the corresponding polynomial functions on $I$, which is a subset of $\Dcal^n(I)$. Since $m_n(f)=n!$ for $f \in \Pcal_n$, one has
\begin{equation}\label{220809a}
D^*(n,p,n!,I)\ioe D^{**}(n,p,I),
\end{equation}
where
\[
D^{**}(n,p,I)=\inf\set{\norm{Q}_p, \; Q\in \Pcal_n}.
\]

\smallskip

A basic fact in the study of the extremal problem $\Ecal^*(n,p,\lambda, I)$ is that \eqref{220809a} is in fact an equality. 

\begin{prop}\label{220809b}
For all $n,p,I$, one has $D^*(n,p,n!,I) = D^{**}(n,p,I)$.
\end{prop}

It follows from this proposition that $C^*(n,p,I)=n!/D^{**}(n,p,I)$ and, by \eqref{220817b},
\begin{equation}\label{220817c}
C(n,p)=L^{n+1/p}n!/D^{**}(n,p,I).
\end{equation}

\smallskip

Let us review the history of Proposition \ref{220809b}.

For $p=\infty$, it is a corollary to a theorem of S. N. Bernstein from 1937. Denoting by $E_k(f)$ the distance (for the uniform norm on $I$) between $f$ and the set of polynomials of degree at most~$k$, he proved in particular that
\[
E_{n-1}(f_0)>E_{n-1}(f_1) \quad (f_0,f_1 \in \Dcal^n(I)),
\]
provided that the inequality $f_0^{(n)}(\xi) > \vert f_1^{(n)}(\xi)\vert$ is valid for every $\xi \in \mathring{I}$ (cf. \cite{bernstein1937}, p. 48, inequalities~(47bis)-(48bis)). Proposition \ref{220809b} follows by taking $f_1(x)=x^n$ and $f_0(x)=\lambda f(x)$, where $f$ is a generic element of~$\Dcal^n(I)$ such that $f^{(n)}(t) \soe n!$ for $a<t<b$, and $\lambda >1$, then letting $\lambda \vers 1$.

This theorem of Bernstein was generalized by Tsenov in 1951 to the case of the $L^p$-norm on~$I$, where $p\soe 1$ (cf. \cite{zbMATH03064306}, Theorem 4, p. 477), thus providing a proof of Proposition \ref{220809b} for $p\soe 1$. The case $0<p<1$ was left open by Tsenov.

The study of the extremal problem $\Ecal^*(n,p,\lambda, I)$ was one of the themes of the 1993 PhD thesis of Xiaoming Huang \cite{MR2689365}. In Lemma 2.0.7, pp. 9-10, she gave another proof (due to Saff) of Proposition~\ref{220809b} in the case $p=\infty$. For $1\ioe p < \infty$, she gave a proof of Proposition 1 which is unfortunately incomplete (cf. \cite{MR2689365}, pp. 28-30). Again, the case $0<p<1$ was left open.

\medskip

We present now a self-contained proof of Proposition~\ref{220809b}, valid for $0<p\ioe \infty$. As it proceeds by induction on $n$, we will need the following classical-looking division lemma, for which we could not locate a reference (compare with \cite{zbMATH03108045} or \cite{zbMATH03316213}).

\begin{prop}\label{220811a}
Let $n\soe 2$ and $f \in \Dcal^n(I)$. Let $c \in \, [a,b]$. Put
\begin{equation}
g(x)=
\begin{dcases}\label{220811a}
\frac{f(x)-f(c)}{x-c} & (x \in I, \, x \neq c)\\
f'(c) & (x=c).
\end{dcases}
\end{equation}

Then $g\in \Dcal^{n-1}(I)$. For every $x \in \,]a,b[\,$, one has
\[
g^{(n-1)}(x)=\frac{f^{(n)}(\xi)}n\virg
\]
where $\xi \in \,]a,b[$.
\end{prop}
\dem

Since $f'$ is continuous, one has
\[
g(x)=\int_0^1f'\big(c+t(x-c) \big) \, dt \quad (x \in I).
\]

Using the rule of differentiation under the integration sign, one sees that $g$ is $n-2$ times differentiable on $I$, with
\[
g^{(n-2)}(x)=\int_0^1t^{n-2}f^{(n-1)}\big(c+t(x-c) \big) \, dt \quad (x \in I).
\]

As $f^{(n-1)}$ is continuous on $I$, this formula yields the continuity of $g^{(n-2)}$ on $I$.

\smallskip

The function $g$ is $n$ times differentiable on $\mathring{I} \setminus\set{c}$ (this set is just $\mathring{I}$ if $c=a$ or $c=b$), being a quotient of $n$ times differentiable functions, with non-vanishing denominator. In the case~$a<c<b$, we have now to check that $g$ is  $n-1$ times differentiable at the point $c$.

The function $f^{(n-1)}$ being continuous on $I$ and differentiable at the point $c$, there exists a function $\eps(h)$, defined and continuous on the segment $[a-c,b-c]$ (the interior of which contains~$0$), vanishing for $h=0$, such that
\[
f^{(n-1)}(c+h)=f^{(n-1)}(c)+hf^{(n)}(c)+h\eps(h) \quad (a \ioe c+h \ioe b).
\]

Hence,
\begin{align*}
g^{(n-2)}(x)&=\int_0^1t^{n-2}f^{(n-1)}\big(c+t(x-c) \big) \, dt \\
&= \int_0^1t^{n-2}\Big(f^{(n-1)}(c)+t(x-c)f^{(n)}(c)+t(x-c) \eps\big(t(x-c)\big) \Big)\, dt \\
&=\frac{f^{(n-1)}(c)}{n-1}+\frac{f^{(n)}(c)}{n}(x-c)+(x-c)\int_0^1t^{n-1}\eps\big(t(x-c)\big) \Big)\, dt 
\end{align*}

When $x$ tends to $c$, the last integral tends to $0$, so that the function $g^{(n-2)}$ is differentiable at the point $c$, with
\[
g^{(n-1)}(c)=\frac{f^{(n)}(c)}{n}\cdotp
\]

\smallskip

If $x \in \mathring{I} \setminus\set{c}$, one may use the general Leibniz rule and Taylor's theorem with the Lagrange form of the remainder in order to compute $g^{(n-1)}(x)$ :
\begin{align*}
g^{(n-1)}(x) &= \frac{d^{n-1}}{dx^{n-1}} \Big(\big(f(x)-f(c)\big)\cdot \frac{1}{x-c}\Big)\\
&=\big(f(x)-f(c)\big)\cdot \frac{(-1)^{n-1}(n-1)!}{(x-c)^n}+\sum_{k=1}^{n-1}\binom{n-1}{k}f^{(k)}(x)\cdot \frac{(-1)^{n-1-k}(n-1-k)!}{(x-c)^{n-k}}\\
&=\frac{(n-1)!}{(c-x)^n}\Big(f(c)-f(x)-\sum_{k=1}^{n-1} \frac{f^{(k)}(x)}{k!}(c-x)^k\Big)\\
&=\frac{(n-1)!}{(c-x)^n}\cdot \frac{f^{(n)}(\xi)}{n!}(c-x)^n\expli{where $\xi$ belongs to the open interval bounded by $c$ and $x$}\\
&=\frac{f^{(n)}(\xi)}n\cdotp\fine
\end{align*}

\smallskip

In the next proposition, we stress the main element of our proof of Proposition \ref{220809b}, namely the fact that the condition $f^{(n)} \soe n!$, for some $f \in \Dcal^n(I)$, implies that the absolute value of $f$ dominates the absolute value of some monic polynomial of degree $n$.

\begin{prop}\label{220811b}
Let $n\soe 1$ and $f \in \Dcal^n(I)$ such that $f^{(n)}(x) \soe n!$ for every $x \in \, ]a,b[\,$.

\smallskip

Then there exists a monic polynomial $P$ of degree $n$, with all its zeros in $I$, such that the inequality $\abs{f(x)}\soe \abs{P(x)}$ is valid for every $x \in I$.

\smallskip

Moreover, if $\abs{f(x)} = \abs{Q(x)}$ for every $x \in I$, where $Q$ is a monic polynomial of degree $n$ with real coefficients, then $f(x)=Q(x)$ for every $x \in I$.
\end{prop}
\dem

The assertion about the zeros may be obtained \emph{a posteriori}, by replacing the zeros of $P$ by their projections on $I$. The following proof leads directly to a polynomial $P$ with all zeros in $I$.

We use induction on $n$.

\smallskip

For $n=1$, the function $f$ is continuous on $[a,b]$, differentiable on $]a,b[\,$, with $f'(x)\soe 1$ for~$a<x<b$. 

If $f(a)\soe 0$, one has, for $a<x\ioe b$, $f(x)=f(a)+(x-a)f'(\xi)$ (where $a<\xi <x$), thus~$f(x) \soe x-a$. Hence, one has $\abs{f(x)}\soe \abs{x-a}$ for every $x \in I$.

If $f(b) \ioe 0$, one proves similarly that $\abs{f(x)}\soe \abs{x-b}$ for every $x \in I$.

If $f(a)< 0<f(b)$, there exists $c \in \, ]a,b[$ such that $f(c)=0$. One has then, for every $x \in I$, 
\[
f(x)=f(x)-f(c)=(x-c)f'(\xi) \quad (\text{where }a < \xi< b).
\]
Hence $\abs{f(x)} \soe \abs{x-c}$ for every $x \in I$, and the result is proven for $n=1$. 

\smallskip

Let now $n \soe 2$, and suppose that the result is valid with $n-1$ instead of $n$. Let~$f \in \Dcal^n(I)$ such that $f^{(n)}(x) \soe n!$ for every $x \in \, ]a,b[\,$.

If $f$ vanishes at some point $c\in I$, it follows from Proposition \ref{220809a} that the function $g$ defined on $I$ by
\[
g(x)=
\begin{dcases}
\frac{f(x)}{x-c} & (x \in I, \, x \neq c)\\
f'(c) & (x=c),
\end{dcases}
\]
belongs to $\Dcal^{n-1}(I)$ and that, for every $x \in \,]a,b[\,$, one has
\[
g^{(n-1)}(x)=\frac{f^{(n)}(\xi)}n\virg
\]
where $\xi \in \,]a,b[$, thus  $g^{(n-1)}(x) \soe (n-1)!$. By the induction hypothesis, there exists a monic polynomial~$Q$ of degree $n-1$, with all its roots in $I$, such that $\abs{g(x)}\soe \abs{Q(x)}$ for every $x\in I$. Hence, one has the inequality $\abs{f(x)}\soe \abs{P(x)}$ for every $x \in I$, where $P(x)=(x-c)Q(x)$ is a monic polynomial of degree $n$, with all its roots in $I$.

If $f>0$, it reaches a minimum at some point $c\in I$. Again, it follows from Proposition \ref{220809a} that the function $g$ defined on $I$ by \eqref{220811a} satisfies the required hypothesis for degree $n-1$. Thus there exists a monic polynomial~$Q$ of degree $n-1$, with all its roots in $I$, such that $\abs{g(x)}\soe \abs{Q(x)}$ for every $x\in I$. Hence, one has the inequality 
\[
f(x)-f(c)=\abs{f(x)-f(c)} \soe \abs{P(x)} \quad (x \in I),
\]
where $P(x)=(x-c)Q(x)$. It follows that
\[
\abs{f(x)}=f(x)\soe f(c) + \abs{P(x)} > \abs{P(x)} \quad (x \in I)
\]

If $f<0$, the reasoning is similar by considering a point $c\in I$ where $f$ reaches a maximum.

\smallskip

Let us prove the last assertion. The hypothesis $\abs{f}=\abs{P}$ is equivalent to the equality~$f^2=P^2$, that is $(f-P)(f+P)=~0$. The set $E=\set{x\in I, \; f(x)+P(x)=0}$ has empty interior, since~$f^{(n)}(x)+P^{(n)}(x)=0$ on every open subinterval of $E$, whereas $f^{(n)}(x)+P^{(n)}(x)\soe 2n!$ on~$\mathring{I}$. The set $I\setminus E$ is therefore dense in $I$ ; its elements $x$ all verify $f(x)=P(x)$, hence $f=P$ on $I$ by continuity.\fin

\medskip

Proposition \ref{220809b} is an immediate corollary of Proposition \ref{220811b}: by taking $f$ and $P$ as stated there, one has $\abs{f(x)}\soe \abs{P(x)}$ for every $x \in I$, so that
\begin{equation}\label{220811c}
\int_a^b \abs{f(x)}^p \, dx \soe \int_a^b \abs{P(x)}^p \, dx,
\end{equation}
for every $p>0$ (for $p=\infty$: $\max\abs{f} \soe \max\abs{P}$).  

Moreover, if $p<\infty$, equality in \eqref{220811c} implies that $\abs{f}= \abs{P}$ on $I$, hence $f=P$. 

In other words, if $0<p<\infty$, the extremal problem $\Ecal^*(n,p,n!, I)$ has exactly the same solutions (value of the infimum and extremal functions) as the problem $\Ecal^{**}(n,p,I)$ obtained by considering only monic polynomials of degree $n$, which one may even take with all their roots in~$I$. 

For $p=\infty$, our reasoning does not prove that an extremal function for $\Ecal^*(n,p,n!, I)$ (if it exists) must be a polynomial. This is true anyway, as proved by Huang in \cite{MR2689365}, pp. 10-13.

\section{Extremal polynomials}

One may now use the results of the well developed theory of the extremal problem $\Ecal^{**}(n,p,I)$ for polynomials. Thus, since the integral
\[
\int_a^b \abs{(x-x_1)\cdots(x-x_n)}^p \, dx \quad (x_1,\dots,x_n \in I)
\]
(or the value $\max_{x\in I}\abs{(x-x_1)\cdots(x-x_n)}$) is a continuous function of $(x_1,\dots,x_n)$, the compactness of $I^n$ yields the existence of an extremal (polynomial) function for $\Ecal^{**}(n,p,I)$, hence for $\Ecal^{*}(n,p,n!, I)$.

It is a known fact that the polynomial extremal problem $\Ecal^{**}(n,p,I)$ has a unique solution for all~$p \in\, ]0,\infty]$, but there is no proof valid uniformly for all values of $p$. 

$\bullet$ For $p=\infty$, uniqueness was proved by Young in 1907 (cf. \cite{zbMATH02644168}, Theorem 5, p. 340)) and follows from the general theory of uniform approximation (cf. \cite{zbMATH03770219}, Theorem 1.8, p. 28).

$\bullet$ For $1<p<\infty$, as proved by Jackson in 1921 (cf. \cite{zbMATH02601208}, \S6, pp. 121-122), this is a consequence of the strict convexity of the space $L^p(I)$.

$\bullet$ For $p=1$, this is also due to Jackson in 1921 (cf. \cite{zbMATH02601209}, \S 4, pp. 323-326).

$\bullet$ For $0<p<1$, the uniqueness of the extremal polynomial was proved in 1988 by Kro\'o and Saff (cf. \cite{zbMATH04032294}, Theorem 2, p. 184). Their proof uses the uniqueness property for $p=1$ and the implicit function theorem.

\medskip

We will denote by $T_{n,p,I}$ the unique solution of the extremal problem $\Ecal^{**}(n,p,I)$. Uniqueness gives immediately the relation
\[
T_{n,p,I}(a+b-x)=(-1)^nT_{n,p,I}(x) \quad (x \in \Real).
\]

Another property of these polynomials is the fact that all their roots are simple. For $p=1$, this fact was proved by Korkine and Zolotareff in 1873 (cf. \cite{zbMATH02717982}, pp. 339-340), before their explicit determination of the extremal polynomial (see \S\ref{220830a} below), and their proof extends, \emph{mutatis mutandis}, to the case $1<p<\infty$. For $p=\infty$, this is a property of the Chebyshev polynomials of the first kind (see \S\ref{221108a} below). Lastly, for $0<p<1$, this was proved by Kro\'o and Saff in~\cite{zbMATH04032294},~p.~187.

Define $T_{n,p}=T_{n,p,[-1,1]}$, and write $n=2k+\eps$, where $k \in \Nat$ and $\eps\in\set{0,1}$. It follows from the mentioned results that
\begin{equation}\label{220817a}
T_{n,p}(x)=x^{\eps} (x^2-x_{n,1}(p)^2) \cdots (x^2-x_{n,k}(p)^2) \quad (x \in \Real),
\end{equation}
where
\[
0<x_{n,1} (p)< \dots < x_{n,k}(p) \ioe 1.
\]

Kro\'o, Peherstorfer and Saff have conjectured that all the $x_{n,k}$ are increasing functions of $p$ (cf. \cite{zbMATH04023942}, p. 656,  and \cite{zbMATH04032294}, p. 192).

\section{Results on $C(n,p)$}

\subsection{The case $n=1$}

The value $n=1$ is the only one for which $C(n,p)$ is explicitly known for all $p$.

\begin{prop}\label{220801a}
One has $C(1,p)=2(p+1)^{1/p}$ for $0<p<\infty$, and $C(1,\infty)=2$.
\end{prop}
\dem

By \eqref{220817a}, one has $T_{1,p}(x)=x$, so that, for $0<p<\infty$,
\[
D^{**}(1,p,[-1,1])=\Big(\int_{-1}^1\abs{t}^p\, dt \Big)^{1/p} =\big(2/(p+1)\big)^{1/p},
\]
and, by \eqref{220817c},
\begin{equation*}
C(1,p)=2^{1+1/p}/D^{**}(1,p,[-1,1])=2(p+1)^{1/p}.\fine
\end{equation*}

\smallskip

Note that the Lemma 1.1, p. 6 of \cite{zbMATH00205893}, asserts that $C(1,p) \ioe 2\cdot 3^{1/p}$ for $p \soe 2$, and that bound is~$< 2(p+1)^{1/p}$ for $p>2$. 

\subsection{The case $p=\infty$}\label{221108a}

This is the classical case, solved by Chebyshev in 1853 by introducing the polynomials $T_n$ defined by the relation $T_n(\cos t)=\cos nt$ (now called Chebyshev polynomial of the first kind): the unique solution of the extremal problem $\Ecal^{**}(n,\infty,[-1,1])$ is $2^{1-n}T_n$. Let us record a short proof of this fact.
 
 \smallskip

Take $I=[-1,1]$ and suppose that $P$ is a monic polynomial of degree $n$ satisfying the inequality $\norm{P}_{\infty}\ioe \norm{2^{1-n}T_n}_{\infty}=2^{1-n}$. Then, for $\lambda>1$ the polynomial
\[
Q_{\lambda}=\lambda 2^{1-n}T_n-P
\]
is of degree $n$, with leading coefficient~$\lambda -1$. Moreover, it satisfies 
\[
(-1)^kQ_{\lambda}(\cos k\pi/n)=\lambda 2^{1-n} -(-1)^kP(\cos k\pi/n)>0 \quad (k=0,\dots,n)
\]
By the intermediate value property, $Q_{\lambda}$ has at least $n$ distinct roots, hence exactly $n$, and these roots, say $x_1,\dots, x_n$, have absolute value not larger than $1$. Hence,
\[
\abs{Q_{\lambda}(x)}= (\lambda-1)\abs{(x-x_1)\cdots(x-x_n)}\ioe (\lambda-1)(1+\abs{x})^n \quad (x \in \Real).
\]
When $\lambda \vers 1$, $Q_{\lambda}(x)$ tends to $0$ for every real $x$, which means that $P=2^{1-n}T_n$.

\medskip

One deduces from this theorem the value of $C(n,\infty)$. One has
\[
D^{**}(n,\infty,[-1,1])=\max_{\abs{x} \ioe 1}\abs{2^{1-n}T_n(x)}=2^{1-n},
\]
hence
\begin{equation}\label{220830b}
C(n,\infty)=2^n\cdot n!/D^{**}(n,\infty,[-1,1])=2^{2n-1}n!
\end{equation}
(compare with the upper bound $C(n,\infty) \ioe 2^{n(n+1)/2}n^n$ of \cite{zbMATH03280851}, 3 (a), p. 185). This result is essentially due to Bernstein (1912, cf. \cite{zbMATH02616863}, p. 65).

\smallskip

Qualitatively, the result expressed by \eqref{220830b} was nicely described by Soula in \cite{zbMATH03006312}, p. 86, as follows.
\begin{quote}
Bernstein's principle: the minimum of the absolute value of the $n$-th derivative of an $n$ times differentiable function and the maximum of the absolute value of the $n$-th derivative of an analytic function have similar orders of magnitude.
\end{quote}

\subsection{The case $p=2$}

In this case, the extremal problem $\Ecal^{**}(n,2,[-1,1])$ is an instance of the general problem of computing the orthogonal projection of an element of a Hilbert space onto a finite dimensional subspace. Here, the Hilbert space is $L^2(-1,1)$, the element is the monomial function $x^n$, and the subspace is the set of polynomial functions of degree less than $n$. The solution follows from the theory of orthogonal polynomials: the extremal polynomial for $\Ecal^{**}(n,2,[-1,1])$ is
\[
\frac{2^n(n!)^2}{(2n)!}P_n(x) \quad (\abs{x}\ioe 1),
\]
where $P_n$ is the $n$-th Legendre polynomial, defined by
\[
P_n(x)=\frac{1}{2^nn!}\frac{d^n}{dx^n}(x^2-1)^{n}.
\]

Hence,
\[
D^{**}(n,2,[-1,1])=\frac{2^n(n!)^2}{(2n)!}\norm{P_n}_2=\frac{2^n(n!)^2}{(2n)!}\sqrt{\frac{2}{2n+1}}\virg
\]
(see \cite{zbMATH02581422}, \S15$\cdot$14, p. 305) and
\begin{equation}\label{220830c}
C(n,2)=2^{n+\demi}\cdot n!/D^{**}(n,2,[-1,1])=\frac{(2n)!}{n!}\sqrt{2n+1},
\end{equation}
a result given by Soula in 1932 (cf. \cite{zbMATH03006312}, pp. 87-88).

\subsection{The case $p=1$}\label{220830a}

The problem $\Ecal^{**}(n,1,[-1,1])$ was solved by Korkine and Zolotareff in \cite{zbMATH02717982}: the extremal polynomial is $2^{-n}U_n(x)$, where $U_n$ is the $n$-th Chebyshev polynomial of the second kind, defined by the relation $U_n(\cos t)=\sin (n+1)t/\sin t$.

Therefore, one has
\begin{align*}
D^{**}(n,1,[-1,1]) &=2^{-n}\int_{-1}^1\abs{U_n(x)} \, dx
=2^{-n}\int_{0}^{\pi}\abs{U_n(\cos t)} \, \sin t \,dt\\
&=2^{-n}\int_{0}^{\pi}\abs{\sin (n+1)t}\,dt
=2^{-n}\int_{0}^{\pi}\sin u\,du\\
&=2^{1-n},
\end{align*}
and
\begin{equation}\label{220830d}
C(n,1)=2^{n+1}\cdot n!/D^{**}(n,1,[-1,1])=2^{2n} n!.
\end{equation}

\subsection{Bounds for $C(n,p)$}

We begin with a simple monotony result.

\begin{prop}
For every positive integer $n$, the function $p \mapsto C(n,p)$ is decreasing on the interval $0<p\ioe \infty$.
\end{prop}
\dem

Let $I=[0,1]$. Equivalently, we will see that the function $p \mapsto D^{**}(n,p,I)$ is increasing. This is due to the fact that, for a fixed $f \in L^{\infty}(I)$ such that $\abs{f}$ is not equal almost everywhere to a constant, the function $p \mapsto \norm{f}_p$ is increasing (a consequence of Hölder's inequality). Thus, for every $Q \in \Pcal_n$ and $0<p < p'\ioe \infty$,
\[
 \norm{Q}_{p'} > \norm{Q}_p \soe D^{**}(n,p,I),
 \]
 which implies that $D^{**}(n,p',I) > D^{**}(n,p,I)$.\fin
 
 \smallskip
 
 In particular, \eqref{220830b} and \eqref{220830d} yield the inequalities
 \[
 2^{2n-1}n! < C(n,p) < 2^{2n}n! \quad (1 < p < \infty).
 \]
 
\smallskip

The next proposition implies that the limit of $C(n,p)$ when $p$ tends to $0$ is $(2e)^nn!$.

\begin{prop}
For every positive integer $n$ and every positive real number $p$, one has
\[
2^n(1+np)^{1/p}n!\ioe C(n,p) \ioe (2e)^nn! 
\]
\end{prop}
\dem
 
Equivalently, we will prove that
\begin{equation}\label{220903a}
(2e)^{-n} \ioe D^{**}(n,p,I) \ioe 2^{-n}(1+np)^{-1/p},
\end{equation}
where $I=[0,1]$.

\smallskip

Let $Q(t)=(t-x_1)\cdots(t-x_n)$, where $0 \ioe x_1,\dots, x_n \ioe 1$. One has
\begin{align*}
\ln\norm{Q}_p  &=\frac 1p \ln\int_{0}^1 \abs{Q(t)}^p \, dt\\
&\soe \frac 1p\int_{0}^1 \ln \big(\abs{Q(t)}^p\big) \, dt\expli{by Jensen's inequality}\\
&=\int_{0}^1 \ln \abs{Q(t)} \, dt\\
&=\sum_{k=1}^n\int_{0}^1 \ln \abs{t-x_k} \, dt.
\end{align*}

Now,
\[
\int_{0}^1 \ln \abs{t-x} \, dt=(1-x)\ln(1-x)+x\ln x-1 \quad (0 \ioe x \ioe 1),
\]
attains its minimal value, namely $-1-\ln 2$,  when $x=1/2$. This implies the first inequality of~\eqref{220903a}.

\smallskip
 
To prove the second inequality of~\eqref{220903a}, we just compute $\norm{Q}_p^p $ when $Q(t)=(t-1/2)^n$ :
\begin{equation*}
\int_0^1\abs{t-1/2}^{np}\, dt= 2\frac{(1/2)^{np+1}}{np+1}\cdotp\fine
\end{equation*}

\smallskip

 For $0<p<1$, we can also prove the following result.
 \begin{prop}
 Let $n$ be a positive integer, and $p$ such that $0<p<1$. One has
 \[
1  \ioe \frac{C(n,p)}{ 2^{2n}n! }\ioe \demi(8/\pi)^{1/p}.
\]
\end{prop}
\dem

The first inequality is just $C(n,1) \ioe C(n,p)$. 
 
 To prove the second inequality, let $r$ and $s$ such that $1<s<2$ and $r^{-1}+s^{-1}=1$. Define
 \begin{align*}
 I_1(s) &=\int_{-1}^1 \frac{dt}{(1-t^2)^{s/2}}\\
 I_2(s) &= \int_{-1}^1 \abs{t}^{(s-1)/s} \, \frac{dt}{\sqrt{1-t^2}}\cdotp
 \end{align*}
 
The integrals $I_1(s)$ and $I_2(s)$ may be computed, using the eulerian identity
\[
\int_0^1t^{x-1}(1-t)^{y-1} \, dt= {\mathrm B}(x,y)=\frac{\Gamma(x)\Gamma(y)}{\Gamma(x+y)}\quad (x>0, \; y>0).
\]
The results are
\begin{align*}
I_1(s) &=2^{1-s}\frac{\Gamma(1-\frac s2)^2}{\Gamma(2-s)}\\
I_2(s)&=\frac{\Gamma\big(1-\frac{1}{2s}\big)\Gamma(\demi)}{\Gamma(\frac 32-\frac{1}{2s})}\cdotp
\end{align*}

 Now, let~$Q \in \Pcal_n$ and put~$p'=p/r$. By Hölder's inequality, one has
 \begin{align*}
\int_{-1}^1 \abs{Q(t)}^{p'} \, \frac{dt}{\sqrt{1-t^2}} &\ioe \Big(\int_{-1}^1 \abs{Q(t)}^{p'r} \, dt \Big)^{1/r}\Big(\int_{-1}^1 \frac{dt}{(1-t^2)^{s/2}}\Big)^{1/s}\\
&=\norm{Q}_p^{p'} I_1(s)^{1/s}.
\end{align*}

It was proved by Kro\'o and Saff (cf. \cite{zbMATH04032294}, pp. 182-183) that
\begin{align*}
2^{(n-1)p'}\int_{-1}^1 \abs{Q(t)}^{p'} \, \frac{dt}{\sqrt{1-t^2}}  &\soe \int_{-1}^1 \abs{T_n(t)}^{p'} \, \frac{dt}{\sqrt{1-t^2}} \\
&=\int_0^{\pi} \abs{\cos nu}^{p'} \, du=\int_0^{\pi} \abs{\cos u}^{p'} \, du\\
&=\int_{-1}^1 \abs{t}^{p'} \, \frac{dt}{\sqrt{1-t^2}} \\
&\soe \int_{-1}^1 \abs{t}^{1/r} \, \frac{dt}{\sqrt{1-t^2}} \expli{one has $p'=p/r<1/r$}\\
&=I_2(s).
\end{align*}

Therefore, with $I=[-1,1]$,
\begin{equation}\label{220903b}
\norm{Q}_p \soe 2^{1-n}I_2(s)^{1/p'}I_1(s)^{-1/p's}=2^{1-n}A(s)^{1/p} \quad (1<s<2),
\end{equation}
where
\[
A(s)=I_2(s)^{s/(s-1)}I_1(s)^{-1/(s-1)}.
\]

Hence
\[
A(s)=2\Big(\frac{\Gamma\big(1-\frac{1}{2s}\big)^s\Gamma(\demi)^s\Gamma(2-s)}{\Gamma(1-\frac s2)^2\Gamma\big(\frac 32-\frac{1}{2s}\big)^s}\Big)^{1/(s-1)} \quad (1<s <2).
\]

Putting $f(s)=\ln\Gamma(s)$, one has
\[
\ln A(s)=\ln 2+\frac{sf(1-1/2s)+sf(1/2)+f(2-s)-2f(1-s/2)-sf(3/2-1/2s)}{s-1}\cdotp
\]

When $s$ tends to $1$, the last fraction tends to
\begin{equation*}
\ln \pi +\frac 32 \psi(1/2) -\frac 32 \psi(1)=\ln \pi -3\ln 2,
\end{equation*}
with the usual notation $f'=\Gamma'/\Gamma=\psi$. It follows that
\[
A(s) \vers \frac {\pi}4 \quad (s \vers 1).
\]

Together with \eqref{220903b}, this gives the inequality 
\[
D^{**}(n,p,[-1,1])\soe 2^{1-n}(\pi/4)^{1/p} 
\]
and \eqref{220817c} now implies
\begin{equation*}
C(n,p) \ioe 2^{2n-1}n!(8/\pi)^{1/p}.\fine
\end{equation*}

 \medskip
 
 We now prove an inequality involving three values of the function $C$.
 
 \begin{prop}
 Let $p,q,r$ be positive real numbers such that 
 \[
 \frac 1p =\frac 1q + \frac 1r\cdotp
 \]
 
 Let $m$ and $n$ be positive integers. Then,
\[
\frac{C(m+n,p)}{(m+n)!} \soe \frac{C(m,q)}{m!}\cdot \frac{C(n,r)}{n!}\cdotp
\]
 \end{prop}
\dem
 
Equivalently, by \eqref{220817c}, one has to prove that
\[
D^{**}(m+n,p,I)\ioe D^{**}(m,q,I)\cdot D^{**}(n,r,I),
\]
where $I$ is a segment of the real line.

In fact, if $P\in \Pcal_{m}$ and $Q\in \Pcal_{n}$, then$PQ \in \Pcal_{m+n}$ hence
\[
D^{**}(m+n,p,I)^p\ioe \int_{I} \abs{P(t)Q(t)}^{p} \, dt\ioe \Big(\int_{I} \abs{P(t)}^{q} \, dt \Big)^{p/q} \cdot \Big(\int_{I} \abs{Q(t)}^{r} \, dt \Big)^{p/r} 
\]
by the definition of $D^{**}(m+n,p,I)$ and Hölder's inequality. The greatest lower bound of the last term, when $P$ runs over $\Pcal_m$ and $Q$ runs over $\Pcal_n$, is
\[
D^{**}(m,q,I)^p\cdot D^{**}(n,r,I)^p.
\]

The result follows.\fin

\subsection{An open question}

Finally, observing that
\[
C(n,2) \sim \sqrt{\frac{2}{\pi}}\cdot 2^{2n}n! \quad ( n \vers \infty),
\]
(an exercise on Stirling's formula from \eqref{220830c}), we ask the following question.

\begin{quote}
Is it true that, for every $p>0$, the quantity $2^{-2n}C(n,p)/n!$ tends to a limit when $n$ tends to infinity?
\end{quote}

\bibliographystyle{smfplain}

\bibliography{journal}{}

\providecommand{\bysame}{\leavevmode ---\ }
\providecommand{\og}{``}
\providecommand{\fg}{''}
\providecommand{\smfandname}{et}
\providecommand{\smfedsname}{\'eds.}
\providecommand{\smfedname}{\'ed.}
\providecommand{\smfmastersthesisname}{M\'emoire}
\providecommand{\smfphdthesisname}{Th\`ese}
\begin{thebibliography}{10}

\bibitem{zbMATH02616863}
{\scshape S.~Bernstein} -- {\og Sur l'ordre de la meilleure approximation des
  fonctions continues par des polynomes de degr{\'e} donn{\'e}.\fg},
  \emph{M\'em. Cl. Sci. Acad. Roy. Belg.} \textbf{IV} (1912), no.~1-104.

\bibitem{bernstein1937}
{\scshape S.~N. Bernstein} -- \emph{Extremal properties of polynomials and the
  best approximation of continuous functions of a single real variable. part
  i}, in russian \smfedname, G. R. O . L, Leningrad, Moscow, 1937.

\bibitem{zbMATH02715424}
{\scshape G.~{Darboux}} -- {\og {M\'emoire sur les fonctions
  discontinues.}\fg}, \emph{{Ann. de l'\'Ec. Norm. (2)}} \textbf{4} (1875),
  p.~57--112.

\bibitem{zbMATH03280851}
{\scshape J.~Dieudonn{\'e}} -- \emph{Foundations of modern analysis}, Pure and
  Applied Mathematics, vol.~10, {Academic} {Press}, New {York}-{London}, 1969.

\bibitem{MR2689365}
{\scshape X.~Huang} -- {\og On extremal properties of algebraic
  polynomials\fg}, \smfphdthesisname, The Ohio State University, 1993.

\bibitem{zbMATH02601209}
{\scshape D.~Jackson} -- {\og Note on a class of polynomials of
  approximation.\fg}, \emph{Trans. Am. Math. Soc.} \textbf{22} (1921),
  p.~320--326.

\bibitem{zbMATH02601208}
\bysame , {\og On functions of closest approximation.\fg}, \emph{Trans. Am.
  Math. Soc.} \textbf{22} (1921), p.~117--128.

\bibitem{zbMATH02717982}
{\scshape A.~Korkine {\normalfont \smfandname} G.~Zolotareff} -- {\og Sur un
  certain minimum\fg}, \emph{Nouv. {Ann}.} \textbf{12} (1873), p.~337--356.

\bibitem{zbMATH04023942}
{\scshape A.~Kro{\'o} {\normalfont \smfandname} F.~Peherstorfer} -- {\og On the
  zeros of polynomials of minimal {{\(L_ p\)}}-norm\fg}, \emph{Proc. Am. Math.
  Soc.} \textbf{101} (1987), p.~652--656.

\bibitem{zbMATH04032294}
{\scshape A.~Kro{\'o} {\normalfont \smfandname} E.~B. Saff} -- {\og On
  polynomials of minimal {{\(L_ q\)}}-deviation, {{\(0<q<1\)}}\fg}, \emph{J.
  Lond. Math. Soc., II. Ser.} \textbf{37} (1988), p.~182--192.

\bibitem{zbMATH00205893}
{\scshape M.~K. {Kwong} {\normalfont \smfandname} A.~{Zettl}} -- \emph{{Norm
  inequalities for derivatives and differences}}, {Lecture Notes in
  Mathematics}, vol. 1536, Springer, 1992.

\bibitem{zbMATH03770219}
{\scshape T.~J. Rivlin} -- \emph{An introduction to the approximation of
  functions. {Corr}. reprint of the 1969 orig}, Dover Publications Inc.,
  Mineola, NY, 1981.

\bibitem{zbMATH03316213}
{\scshape L.~Schoenfeld} -- {\og On the differentiability of indeterminate
  quotients\fg}, \emph{Math. Mag.} \textbf{41} (1968), p.~152--155.

\bibitem{zbMATH03006312}
{\scshape J.~Soula} -- {\og Sur une in{\'e}galit{\'e} verifiee par une fonction
  et sa d{\'e}riv{\'e}e d'ordre $n$\fg}, \emph{Mathematica, Cluj} \textbf{6}
  (1932), p.~86--88.

\bibitem{zbMATH03064306}
{\scshape I.~V. Tsenov} -- {\og On a question of the approximation of functions
  by polynomials\fg}, \emph{Mat. Sb., Nov. Ser.} \textbf{28} (1951),
  p.~473--478 (Russian).

\bibitem{zbMATH03108045}
{\scshape H.~Whitney} -- {\og Differentiability of the remainder term in
  {Taylor}'s formula\fg}, \emph{Duke Math. J.} \textbf{10} (1943), p.~153--158.

\bibitem{zbMATH02581422}
{\scshape E.~T. {Whittaker} {\normalfont \smfandname} G.~N. {Watson}} --
  \emph{{A course of modern analysis}}, 4th \smfedname, Cambridge University
  Press, 1927.

\bibitem{zbMATH02644168}
{\scshape J.~W. Young} -- {\og General theory of approximation by functions
  involving a given number of arbitrary parameters.\fg}, \emph{Trans. Am. Math.
  Soc.} \textbf{8} (1907), p.~331--344.

\end{thebibliography}

\medskip

\footnotesize

\noindent BALAZARD, Michel\\
Institut de Mathématiques de Marseille (I2M)\\
CNRS, Aix Marseille Universit\'e \\
Marseille, France\\
e-mail address: \texttt{balazard@math.cnrs.fr}

\end{document}